\documentclass[12pt]{article}
\usepackage{graphicx}
\usepackage{amsmath,amsthm,amssymb,enumerate}
\usepackage{euscript,mathrsfs}
\usepackage{color}
\usepackage{dsfont}
\usepackage{url}
\usepackage[left=2cm,right=2cm,top=3.5cm,bottom=3.5cm]{geometry}
\usepackage{color}
\usepackage[framemethod=tikz]{mdframed}
\allowdisplaybreaks

\usepackage{soul}

\catcode`\@=11 \@addtoreset{equation}{section}

\catcode`\@=12

\newtheorem{Theorem}{Theorem}[section]
\newtheorem{Proposition}[Theorem]{Proposition}
\newtheorem{Lemma}[Theorem]{Lemma}
\newtheorem{Corollary}[Theorem]{Corollary}

\theoremstyle{definition}
\newtheorem{Definition}[Theorem]{Definition}

\newtheorem{Remark}[Theorem]{Remark}

\newcommand{\bTheorem}[1]{
	\begin{Theorem} \label{T#1} }
	\newcommand{\eT}{\end{Theorem}}

\newcommand{\bProposition}[1]{
	\begin{Proposition} \label{P#1}}
	\newcommand{\eP}{\end{Proposition}}

\newcommand{\bLemma}[1]{
	\begin{Lemma} \label{L#1} }
	\newcommand{\eL}{\end{Lemma}}

\newcommand{\bCorollary}[1]{
	\begin{Corollary} \label{C#1} }
	\newcommand{\eC}{\end{Corollary}}

\newcommand{\bRemark}[1]{
	\begin{Remark} \label{R#1} }
	\newcommand{\eR}{\end{Remark}}

\newcommand{\bDefinition}[1]{
	\begin{Definition} \label{D#1} }
	\newcommand{\eD}{\end{Definition}}

\newcommand{\vrB}{\vr_B}

\newcommand{\Ds}{\mathbb{D}_x}

\newcommand{\vuB}{\vc{u}_B}

\newcommand{\bfphi}{\boldsymbol{\varphi}}

\newcommand{\bFormula}[1]{
	\begin{equation} \label{#1}}
	\newcommand{\eF}{\end{equation}}

\newcommand{\Ov}[1]{\overline{#1}}

\newcommand{\aleq}{\stackrel{<}{\sim}}

\newcommand{\ageq}{\stackrel{>}{\sim}}

\newcommand{\vr}{\varrho}

\newcommand{\tvr}{\wtilde \vr}
\newcommand{\tvu}{{\wtilde u}}

\newcommand{\vu}{\vc{u}}
\newcommand{\vm}{\vc{m}}

\newcommand{\vc}[1]{{\bf #1}}

\newcommand{\Div}{{\rm div}_x}
\newcommand{\Grad}{\nabla_x}

\newcommand{\dx}{\,{\rm d} {x}}

\newcommand{\dt}{\,{\rm d} t }

\newcommand{\vU}{\vc{U}}

\newcommand{\intO}[1]{\int_{\Omega} #1 \ \dx}

\newcommand{\D}{{\rm d}}

\newcommand{\br}{ \nonumber \\ }

\def\softd{{\leavevmode\setbox1=\hbox{d}%
		\hbox to 1.05\wd1{d\kern-0.4ex{\char039}\hss}}}
\definecolor{Cgrey}{rgb}{0.85,0.85,0.85}
\definecolor{Cblue}{rgb}{0.50,0.85,0.85}
\definecolor{Cred}{rgb}{1,0,0}
\definecolor{fancy}{rgb}{0.10,0.85,0.10}
\definecolor{amaranth}{rgb}{0.9, 0.17, 0.31}

\newcommand\Cbox[2]{%
	\newbox\contentbox%
	\newbox\bkgdbox%
	\setbox\contentbox\hbox to \hsize{%
		\vtop{
			\kern\columnsep
			\hbox to \hsize{%
				\kern\columnsep%
				\advance\hsize by -2\columnsep%
				\setlength{\textwidth}{\hsize}%
				\vbox{
					\parskip=\baselineskip
					\parindent=0bp
					#2
				}%
				\kern\columnsep%
			}%
			\kern\columnsep%
		}%
	}%
	\setbox\bkgdbox\vbox{
		\color{#1}
		\hrule width  \wd\contentbox %
		height \ht\contentbox %
		depth  \dp\contentbox
		\color{black}
	}%
	\wd\bkgdbox=0bp%
	\vbox{\hbox to \hsize{\box\bkgdbox\box\contentbox}}%
	\vskip\baselineskip%
}

\mdfdefinestyle{MyFrame}{%
	linecolor=black,
	outerlinewidth=1pt,
	roundcorner=5pt,
	innertopmargin=\baselineskip,
	innerbottommargin=\baselineskip,
	innerrightmargin=10pt,
	innerleftmargin=10pt,
	backgroundcolor=white!20!white}

\newcommand{\wtilde}{\widetilde}

\allowdisplaybreaks

\begin{document}


\title{\bf Unconditional stability of radially symmetric steady sates of compressible viscous fluids with inflow/outflow boundary conditions}

\author{Eduard Feireisl
	\thanks{The work of E.F. was partially supported by the
		Czech Sciences Foundation (GA\v CR), Grant Agreement
		24--11034S. The Institute of Mathematics of the Academy of Sciences of
		the Czech Republic is supported by RVO:67985840. This work was partially supported by the Thematic Research Programme, University of Warsaw, Excellence Initiative Research University.} \and Piotr Gwiazda
	\and Agnieszka \'Swierczewska-Gwiazda
	\thanks{The work of P. G. and A. \'S-G. was partially supported by  National Science Centre
		(Poland),  agreement no   2023/51/B/ST1/01546.}}

\date{}

\maketitle

\medskip

\centerline{Institute of Mathematics of the Academy of Sciences of the Czech Republic}

\centerline{\v Zitn\' a 25, CZ-115 67 Praha 1, Czech Republic}

\centerline{feireisl@math.cas.cz}

\medskip

\centerline{Institute of Mathematics of Polish Academy of Sciences}
\centerline{\'Sniadeckich 8, 00-956 Warszawa, Poland}
\centerline{pgwiazda@mimuw.edu.pl}

\medskip

\centerline{Institute of Applied Mathematics and Mechanics, University of Warsaw}

\centerline{Banacha 2, 02-097 Warsaw, Poland}
\centerline{aswiercz@mimuw.edu.pl}

\begin{abstract}
	
We show that certain radially symmetric steady states of compressible viscous fluids in domains with inflow/outflow boundary conditions are unconditionally stable. This means that \emph{any} not necessarily radially symmetric solution of the associated evolutionary problem converges to a single radially symmetric steady state.

\end{abstract}


{\small

\noindent
{\bf 2020 Mathematics Subject Classification:}{
(primary) 35 Q 30;
(secondary) 35 B 40, 35 D 30}

\medbreak

\tableofcontents

}

\section{Introduction}
\label{i}

The problem of global asymptotic stability of steady state solutions to problems in fluid mechanics is both mathematically challenging and physically relevant. Our goal is to consider 
compressible viscous fluids described by means of the 
\emph{barotropic Navier--Stokes (NS) system}:
\begin{align}
	\partial_t \vr + \Div (\vr \vu ) &= 0, \label{p1} \\ 
	\partial_t( \vr \vu) + \Div (\vr \vu \otimes \vu) + \Grad p (\vr) &= 
	\Div \mathbb{S}(\Ds \vu), \ \Ds \vu \equiv \frac{1}{2} \left( \Grad \vu + \Grad^t \vu \right) , 
	\label{p2}
\end{align}
for $(t,x) \in (0,T) \times \Omega$, $\Omega \subset R^3$. The unknowns describing the state of the system are the fluid mass density $\vr = \vr(t,x)$ and the velocity $\vu = \vu(t,x)$. 
The pressure $p = p(\vr)$ is related to  the density by a specific equation of state (EOS), where we suppose
\begin{align} \label{p3}
	p \in C[0, \Ov{\vr}) \cap C^1(0, \Ov{\vr}),\ 
	0 < \Ov{\vr} \ {\leq} \infty,\ p(0) = 0,\ p'(\vr) > 0 \ \mbox{for}\ 
	0 < \vr < \Ov{\vr}. 
\end{align}	
{The case $\Ov{\vr} < \infty$ characterizes the so--called hard sphere pressure law. A detailed description of the 
pressure--density equation of state is presented in Section \ref{EOS}.}	
The viscous stress $\mathbb{S}$ is given by \emph{Newton's rheological law}
\begin{equation} \label{p4} 
	\mathbb{S} = \mathbb{S}(\Ds \vu) = \mu \Big( \Grad \vu + \Grad^t \vu - 
	\frac{2}{3} \Div \vu \mathbb{I} \Big) + \lambda \Div \vu,    	
\end{equation}
with constant viscosity coefficients $\mu > 0$, $\lambda \geq 0$.

Our goal is to establish global stability of the radially symmetric steady states. Accordingly, the underlying spatial domain $\Omega$ is radially symmetric. We focus on two basic 
problems corresponding to the outflow and inflow--outflow boundary conditions. 

\begin{itemize}
\item {\bf Outflow problem on exterior domain.}
We consider
\begin{equation} \label{p5}
	\Omega = R^3 \setminus \Ov{B}_{ \Ov{r}},\ \mbox{where} \ B_{\Ov r} = \{ |x| < { {\Ov{r}} }\}, 
\end{equation}
is an exterior domain, and impose the outflow boundary conditions for the velocity, 
\begin{equation} \label{p6} 
	\vu|_{|x| = \Ov{r}} = u_B \vc{n},\ u_B > 0.
\end{equation}
Here and hereafter, the symbol $\vc{n}$ denotes the \emph{outer} normal vector to $\partial \Omega$.

In addition, the far field behaviour is specified, 
\begin{equation} \label{p7} 
	\vu \to 0,\ \vr \to \vr_\infty  \ \mbox{as}\ |x| \to \infty,\ 
	0 < \vr_\infty < \Ov{\vr}.
\end{equation}

\item {\bf Inflow--outflow problem on a bounded domain.}

We consider 
	\begin{equation} \label{p8}
		\Omega = B_{r+} \setminus \Ov{B}_{r^-},\ 
		0 < r^{-} < r^+,
	\end{equation}	
	with the inflow--outflow boundary conditions 
	\begin{align} 
		\vu|_{|x| = r^-} & = - u_B^- \vc{n},\ u_B^- > 0,  
		{\vr}|_{|x| = r^-} = \vr_B > 0, \label{p9} \\ 
		\vu|_{|x| = r^+} &=  u_B^+ \vc{n},\ u_B^+ > 0. \label{p10} 
	\end{align}
	\end{itemize}
{In both cases we also consider a vector field $\vuB$ defined on the whole 
	physical space $R^3$ that coincides with the boundary data on $\partial \Omega$.}

The problem of existence and qualitative behaviour 
of steady states to the NS system with the aforementioned boundary conditions has been extensively studied in a series of papers by Hashimoto and Matsumura \cite{HasMat2}, \cite{HasMat1}, Hashimoto, Nishibata, and Sugizaki \cite{HaNiSu}, or, more recently, Huang et al. \cite{HuHaNi}, \cite{HuNi}. Besides a detailed description of the steady state, their stability with respect to \emph{radially symmetric} perturbations is established. Our goal is to show that stability in fact holds in a much larger class of weak solutions to the NS system, not necessarily radially symmetric. {In contrast with the radially symmetric 
case that can be reduced to a relatively simple 1-D problem, the general 3-D case is much more complex, in particular a suitable existence theory 
in the class of global--in--time large data solutions is necessary.}	 

Our approach is based on a recently developed theory of 
weak solutions to the NS system with general inflow--outflow boundary conditions \cite{FeiNovOpen}, see also 
Choe, Novotn\'y and Yang \cite{ChoNoYa}, Ne\v casov\'a, Novotn\'y and Roy \cite{NeNoRo}. The distance between a solution of the evolutionary problem and the steady state is measured by the 
Bregman divergence (relative energy) {that can be seen as a linear perturbation of the total energy of the system}. In order to establish convergence, the 
coercive terms in the relative energy inequality must be controlled by dissipation. To this end, qualitative properties of the steady state play a crucial role. The relevant piece of information in the outflow exterior domain problem is provided by the results Hashimoto and Matsumura \cite{HasMat2}, Huang and Nishibata \cite{HuNi}. The inflow--outflow problem is considered as a perturbation of the same system of equations and boundary conditions posed on a ``flat'' strip like domain. As an intermediate step, we therefore obtain an unconditional convergence result on a periodic strip. Finally, convergence for the inflow--outflow problem \eqref{p8}--\eqref{p10} results from 
careful analysis of the perturbed steady solution.

The paper is organized as follows. In Section \ref{m}, we recall the concept of weak solution to the NS system, and formulate our main results. Section \ref{s} is devoted to the steady solutions and their basic properties. In Section \ref{r}, we introduce the relative energy, together with the associated relative energy inequality. The proof of main results is completed in Section \ref{c}.

\section{Problem formulation, weak solutions, and main results}
\label{m}

We suppose that $\Omega \subset R^3$ is a domain with compact Lipschitz boundary. The boundary velocity is determined by a 
(smooth) function $\vu_B \in C^1_c(R^3; R^3)$. Similarly, the 
boundary density on the inflow part of the boundary is the trace of $\vr_B \in C^1(R^3)$. We consider the problem on the infinite space--time cylinder 
{$[0,\infty) \times \Omega$}. 

\subsection{Weak formulation}
\label{wf}

The weak formulation of equation of continuity \eqref{p1}, together with the inflow/outflow boundary conditions reads:
\begin{align}
\int_0^\tau	&\intO{ \Big( \vr \partial_t \varphi + \vr \vu \cdot 
	\Grad \varphi \Big) } \dt = \left[ \intO{ \vr \varphi } \right]_{t = 0}^{t = \tau} \br
	&+ \int_0^\tau \int_{\partial \Omega} \varphi \vr_b 
	[\vu_B \cdot \vc{n}]^- \D S_x \dt + \int_0^\tau \int_{\partial \Omega} \varphi \vr [\vu_B \cdot \vc{n}]^+ \D S_x \dt
	\label{m1}
\end{align}
for any test function $\varphi \in C^1_c([0, \infty) \times R^3)$, and {any $\tau > 0$}.
 {Here and hereafter,  
	\[
	[a]^+ = \max\{ 0, a \},\ [a]^- = \min\{0,a\}
	\]
	denote the positive and negative part of $e \in R$, respectively.}
In addition, {we consider the weak solutions} satisfying the renormalized version of \eqref{m1}: 
\begin{align}
	\int_0^\tau	&\intO{ \Big( b(\vr) \partial_t \varphi + b(\vr) \vu \cdot 
		\Grad \varphi + \left( b(\vr) - b'(\vr) \vr \right) \Div \vu \varphi \Big) } \dt = \left[ \intO{ b(\vr) \varphi } \right]_{t = 0}^{t = \tau} \br
	&+ \int_0^\tau \int_{\partial \Omega} \varphi b(\vr_B) 
	[\vu_B \cdot \vc{n}]^- \D S_x \dt + \int_0^\tau \int_{\partial \Omega} \varphi b(\vr) [\vu_B \cdot \vc{n}]^+ \D S_x \dt
	\label{m2}
\end{align}
for any test function $\varphi \in C^1_c([0, \infty) \times R^3)$, any $b \in C^1[0, \infty)$, $b' \in C_c[0, \infty)$, and 
any {$\tau > 0$}.
All quantities are supposed to be {at least} integrable. We refer to 
\cite[Chapter 3, Section 3.3.1]{FeiNovOpen} for the interpretation of the ``trace'' of the function $\vr$ on the outflow component of the boundary. 

The momentum equation \eqref{p2} is satisfied in the sense of distributions:
\begin{align}
\int_0^\tau &\intO{ 
\Big( \vr \vu \cdot \partial_t \bfphi + (\vr \vu \otimes \vu): \Grad \bfphi + p(\vr) \Div \bfphi \Big) } \dt = 
\int_0^\tau	\intO{ \mathbb{S}(\Ds \vu) : \Grad \bfphi } \dt \br
&+ \left[ \intO{ \vr \vu \cdot \bfphi } \right]_{t=0}^{t = \tau}
\label{m3}
\end{align}
for any {$\bfphi \in C^1_c([0, \infty; R^3) \times \Omega)$, $\tau > 0$.} 
{Formula \eqref{m3} requires integrability of the velocity gradient. In particular,
the velocity $\vu$ admits a trace on $\partial \Omega$, and we impose}
\begin{equation} \label{m4}
\vu = \vuB \ \mbox{on}\ \partial \Omega.
\end{equation}
{Similarly to the above}, we tacitly assume that all integrands are well defined and at least locally integrable. In particular, 
\begin{equation} \label{dens}
0 \leq \vr {< \Ov{\vr}} \ \ \mbox{a.e. in}\ (0, \infty) \times \Omega,
\end{equation}	
where {$\Ov{\vr} \leq  \infty$} is specified by the EOS {in 
Section \ref{EOS} below.} {The boundary condition \eqref{m4} 
as well as the far the relevant far field conditions for the velocity $\vu$ are interpreted in the weak sense, 
\[
\vu - \vuB \in L^2_{\rm loc}(0,T; D^{1,2}_0 (\Omega; R^3)), 
\]
where the $D^{1,2}_0(\Omega; R^3)$ is the completion of $C^\infty_c(\Omega; R^3)$ with respect to the gradient norm $\| \Grad \vc{v} \|_{L^2(\Omega; R^{3 \times 3})}$. }

Finally, we consider only the weak solutions satisfying some form of energy balance. The latter can be written in the form:
\begin{align} 
	& \left[ \intO{ E \left(\vr, \vu \Big| \vr_\infty, \vuB \right) } \right]_{t = 0}^{t = \tau} 
	+ \int_0^\tau \intO{ {\mathbb{S} \Big(\Ds \vu - \Ds \vuB \Big): \Big( \Ds \vu - \Ds \vuB \Big) } } \dt \br 
	&+ \int_0^\tau \int_{\partial \Omega} \Big( P(\vr) - P'(\vr_\infty) (\vr - \vr_\infty) - P(\vr_\infty) \Big) [\vuB \cdot \vc{n}]^+ \D \sigma_x \dt	\br 
	&+ \int_0^\tau \int_{\partial \Omega} \Big( P(\vr_B) - P'(\vr_\infty) (\vr_B - \vr_\infty) - P(\vr_\infty) \Big) [\vuB \cdot \vc{n}]^- \D \sigma_x \dt \br
	&\leq - \int_0^\tau \intO{ \vr (\vu - \vuB) \cdot \Grad \vuB \cdot (\vu - \vuB) }
	\br &\quad -   \int_0^\tau \intO{ \Big( p(\vr) - p'(\vr_\infty) (\vr - \vr_\infty) - p(\vr_\infty) \Big) \Div \vuB } \dt \br 
	&\quad + \int_0^\tau \intO{ \vr (\vuB - \vu) \cdot \Div ( \vuB \otimes \vuB)    \Big] } \dt \br 
	&\quad + \int_0^\tau \intO{ \Big( \vr (\vu - \vuB) \cdot \vuB + 
		p'(\vr_\infty) \left( \vr_\infty - \vr \right) \Big)  \Div  \vuB  } \dt \br 
	&\quad + \int_0^\tau \intO{  
		\Div \mathbb{S} {(\Ds \vuB) \cdot} (\vuB - \vu) } \dt	
	\label{m5}	
\end{align}
where we have set 
\[
E \left( \vr, \vu \Big| \vr_\infty, \vu_B \right) = 
\frac{1}{2} |\vu - \vu_B |^2 + P(\vr) - P'(\vr_\infty) (\vr - \vr_\infty) 
- P(\vr_\infty), 
\]
with the pressure potential 
\begin{equation} \label{r2}
P'(\vr) \vr - P(\vr) = p(\vr).
\end{equation}
Note that \eqref{m5} is nothing other than a particular form of a general relative energy inequality introduced in Section \ref{r} below. The constant 
$\vr_\infty$ coincides with the far field value of the density prescribed in \eqref{p7}, and can be taken arbitrary non--negative number if $\Omega$ is a bounded domain.

If $\Omega$ is bounded, it is possible to establish the existence of weak solutions satisfying \eqref{m5} in the differential form \cite[Chapter 3, Section 3.2]{FeiNovOpen}:
\begin{align} 
	\frac{\D }{\dt}&\intO{ E \left(\vr, \vu \Big| \vr_\infty, \vuB \right) }  
	+ \intO{ \mathbb{S} {\Big(\Ds \vu - \Ds \vuB \Big): \Big( \Ds \vu - \Ds \vuB \Big) } }  \br 
	&+ \int_{\partial \Omega} \Big( P(\vr) - P'(\vr_\infty) (\vr - \vr_\infty) - P(\vr_\infty) \Big) [\vuB \cdot \vc{n}]^+ \D \sigma_x 	\br 
	&+ \int_{\partial \Omega} \Big( P(\vr_B) - P'(\vr_\infty) (\vr_B - \vr_\infty) - P(\vr_\infty) \Big) [\vuB \cdot \vc{n}]^- \D \sigma_x \dt \br
	&\leq - \int_0^\tau \intO{ \vr (\vu - \vuB) \cdot \Grad \vuB \cdot (\vu - \vuB) }
	\br &\quad - \intO{ \Big( p(\vr) - p'(\vr_\infty) (\vr - \vr_\infty) - p(\vr_\infty) \Big) \Div \vuB }  \br 
	&\quad +  \intO{ \vr (\vuB - \vu) \cdot \Div ( \vuB \otimes \vuB)    \Big] }  \br 
	&\quad +  \intO{ \Big( \vr (\vu - \vuB) \cdot \vuB + 
		p'(\vr_\infty) \left( \vr_\infty - \vr \right) \Big)  \Div  \vuB  }  \br 
	&\quad + \intO{  \Div \mathbb{S} {(\Ds \vuB) \cdot} (\vuB - \vu) }	
	\label{m6}	
\end{align}
satisfied in $\mathcal{D}'(0, \infty)$. Validity of the differential form {of the energy inequality} 
\eqref{m6} can be extended to unbounded (exterior) spatial domains as long as 
$\vr_\infty = 0$, more exactly {on condition that} the total mass of the fluid is finite,
\[
\intO{ \vr } < \infty,
\]
see \cite{DoroFe}. For $\vr_\infty > 0$ considered in 
\eqref{p7}, the existence of weak solutions satisfying \eqref{m6} is an open problem, cf. also Novotn\' y and Pokorn\' y \cite{NovPok07}.

\subsection{Pressure--density equation of state}
\label{EOS}

The specific form of the equation of state (EOS) relating the pressure $p$ to the density of a barotropic fluid plays a crucial role in the existence theory as well as in the analysis of the long--time behaviour. In addition to 
\eqref{p3}, we suppose 
\begin{equation} \label{m7}
	p \in C^2(0, \Ov{\vr}),\ p''(\vr) > 0 \ \mbox{in}\ (0, \Ov{\vr}),\ 
\lim_{\vr \to \Ov{\vr} -} p(\vr) = \infty.
\end{equation}	

There are two iconic examples of EOS satisfying \eqref{p3}, \eqref{m7}:
\begin{itemize}
	\item {\bf Isentropic EOS.} 
\begin{equation} \label{m8}
p (\vr) = a \vr^\gamma ,\ a > 0, \gamma > 1, \ \Ov{\vr} = \infty.	
\end{equation}
\item {\bf Hard--sphere EOS.}
\begin{equation} \label{m9}
0 < \Ov{\vr} < \infty,\ p(\vr) \approx (\Ov{\vr} - \vr)^{-\beta} 
\ \mbox{for}\ \vr \to \Ov{\vr},\ \beta > 0.
\end{equation}
\end{itemize} 
{Here and hereafter, the symbol $a \approx b$ denotes that  $a \aleq b$ and $b \aleq a$, where
$a \aleq b$ means $a \leq c b$, where $c > 0$ is a generic constant. }

The \emph{existence} of weak solutions of the NS system 
with the isentropic EOS
for $\gamma \geq \frac{9}{5}$ if $d = 3$, and $\gamma \geq \frac{3}{2}$ if $d=2$ was established in the seminal work of Lions \cite{LI4}. The results 
have been later extended to $\gamma > \frac{d}{2}$ in \cite{FNP}. 
As for the hard--sphere EOS, we refer to \cite{FeiZha} and the monograph \cite[Chapter 9]{FeiNovOpen} for the existence theory for $\beta > 3$. 
The inflow--outflow problem for the hard--sphere EOS has been elaborated in detail in Choe, Novotn\' y, and Yang \cite{ChoNoYa}. As they pointed out, 
\emph{global--in--time} weak solutions exist as long as $\beta > \frac{5}{2}$ 
and the boundary velocity $\vuB$ admits an extension 
\begin{equation} \label{m10}
\vuB \in C^1(\Ov{\Omega},R^3),\ \Div \vuB \geq 0, \Div \vuB \not \equiv 0. 
\end{equation}
As we shall see below, condition \eqref{m10} is satisfied in all cases considered in the present paper. Finally, we mention the recent result 
of Ne\v casov\' a, Novotn\' y, and Roy \cite{NeNoRo} concerning the hard--sphere pressure EOS and exterior spatial domains. 

\subsection{Steady solutions}
\label{s}

Steady solutions $(\vr_s, \vu_s)$ satisfy the stationary NS system:
\begin{align}
\Div (\vr_s \vu_s ) &= 0, \label{S1} \\
\Div (\vr_s \vu_s \otimes \vu_s) + \Grad p(\vr_s) &= {\Div \mathbb{S}(\Ds \vu_s) } \label{S2}
\end{align}
in $\Omega$, along with the relevant boundary conditions. In this paper, we focus on strong solutions satisfying \eqref{S1}, \eqref{S2} in the classical sense. 

\subsubsection{Exterior outflow problem}
\label{ss1}

The existence of a radially symmetric steady solution to the NS system 
\eqref{S1}, \eqref{S2} 
on the exterior domain \eqref{p5}, supplemented with the outflow boundary conditions \eqref{p6}, and the far field conditions \eqref{p7} was
established 
by Hashimoto, Nishibata, and Sugizaki \cite{HaNiSu}. More specifically, 
given 
\begin{equation} \label{s1a}
0 < \vr_\infty < \Ov{\vr}, 
\end{equation}
and a sufficiently small outflow velocity 
\begin{equation} \label{s1}
	0 < u_B \leq \Ov{u} \Big( \mu, { \lambda, \Ov{r}},\ \vr_\infty, p'(\vr_\infty) \Big), 
\end{equation}	
there exists a unique radially symmetric solution of the problem {in the form}
\[
\vr_s(x) = \tvr(|x|),\ 
\vu_s(x) = \frac{x}{|x|} \tvu (|x|) ,\ |x| \geq { \Ov{r}}.
\]
{The functions $\tvr = \tvr(r)$, $\tvu = \tvu(r)$ satisfy the equations} 
\begin{align}
\partial_r\left( r^{2} \tilde{\rho} \tilde{u} \right)&= 0,\\
\tilde{\rho} \tilde{u} \partial_r\tilde{u} + \partial_rp(\tilde{\rho}) &= \left( \frac{4}{3} \mu + \lambda \right)  \partial_r\left(\frac{\partial_r( r^{2} \tilde{u} )}{r^{2}}\right), 
\end{align}
with the far-field condition
\begin{equation}
\lim_{r \to \infty} (\tilde{\rho}, \tilde{u})(r) = (\rho_\infty, 0), 
\end{equation}
and the boundary condition
\begin{equation}
\tilde{u}(\Ov{r}) = u_B.
\end{equation}
The solution enjoys the following properties (cf. \cite[Lemma 1.1]{HaNiSu}):
\begin{align}
0 < \tvr &< \vr_\infty, \ 
\Big( \vr_\infty - \tvr (|x|) \Big) \approx u_B^2 |x|^{-4}, 
\partial_r \tvr > 0,\ \partial_r \tvr(|x|) \approx u_B^2 |x|^{-5},	\label{s2}\\
\tvu( \Ov{r}) &= - u_B , \ \tvu < 0,\ \lim_{|x| \to \infty} \tvu(|x|) = 0,\ \partial_r \tvu > 0,\ 
\partial_r \tvu (|x|) \approx u_B |x|^{-3}, { \tvu (|x|) \approx u_B |x|^{-2}}.
\label{s3}
	\end{align}
Note carefully 
\[
0 < \vr_s (|x|) \leq \vr_\infty \ \mbox{for all}\ |x| > {\Ov{r}}. 
\]
In particular, the construction {elaborated in \cite{HaNiSu}} relies only on the local properties of the 
pressure EOS around the far field value $\vr_\infty$ and applies to the hard sphere pressure as long as \eqref{s1a} holds.	

Finally, a routine manipulation yields 
\begin{equation} \label{s4}
\Grad \vu_s = \left( \partial_r \tvu (|x|) - \frac{\tvu(|x|)}{|x|} \right) \frac{ x \otimes x }{|x|^2} + \frac{\tvu(|x|)}{|x|} {\rm Id} .
	\end{equation}
Moreover, exploiting the fact that $\vr_s$, $\vu_s$ solve the equation of continuity, we deduce 
\begin{equation} \label{s5}
\Div \vu_s (x) = - \frac{1}{\vr_s (x) } \Grad \vr_s (x) \cdot \vu_s = 
- \frac{1}{\tvr(|x|)} \partial_r \tvr (|x|) \ \tvu (|x|) > 0. 
\end{equation}
In particular, the condition \eqref{m10}, necessary for the existence of global--in--time solutions to the evolutionary problem, is satisfied for 
$\vuB = \vu_s$.

\subsubsection{Inflow--outflow problem on a periodic strip}
\label{ss3}

Similarly to the preceding section, we look for solutions of the inflow--outflow problem \eqref{p8}--\eqref{p10} in the class of radially symmetric functions. As an intermediate step, we first introduce a similar problem in the simpler geometry of a (periodic) strip. 

Specifically, consider 
\begin{equation} \label{s12}
	\Omega = \left\{ x = (y,r) \ \Big| \ y \in \mathbb{T}^2,\ 
	r \in (0,1) \right\},
\end{equation}	
{where $\mathbb{T}^2$ denotes a 2-D (flat) torus},
with the inflow--outflow conditions 
\begin{align}
	\vu|_{r = 0} &= (0,0, u_B^-),\ \vr|_{r = 0} = \vr_B , \label{s13} \\  
	\vu|_{r = 1} &=(0,0,u_B^+). \label{s14}
\end{align}
The corresponding steady state depends only on the vertical coordinate $r$,
\[
\vr_s = \tvr(r),\ \vu_s = (0,0, \tvu(r)). 
\]
The system of equations reads 
\begin{align}
	\partial_r (\tvr \tvu) &= 0 \ \Rightarrow\ 
	\tvr(r) \tvu(r) = \vr_B u^-_B, \label{s15} \\
	\vr_B \vu^-_B \partial_r \tvu + \partial_r p(\tvr) &= 
	\left( \frac{4}{3} \mu + \lambda \right) \partial^2_{r,r} \tvu. \label{s16}	
\end{align}
Equation \eqref{s16} can be integrated to obtain 
\begin{equation} \label{s17}
	\left( \frac{4}{3} \mu + { \lambda}\right) \partial_r \tvu = 
	p\left( \frac{\vr_B u^-_B}{\tvu}    \right) + \vr_B \vu^-_B \tvu + \Lambda,\ 
	\tvu(0) = u^-_B > 0,
\end{equation}
where the value of the integration constant $\Lambda$ is determined by the outflow velocity
\begin{equation} \label{s18}
	\tvu(1) = u^+_B.
\end{equation}

The first observation is that the unique solution $\tvu$ of \eqref{s17}
remains strictly positive as long as the pressure satisfies 
\begin{equation} \label{s19}
{	\lim_{\vr \to \Ov{\vr}-} p(\vr) = \infty. }
\end{equation}	
This is definitely true for the isentropic EOS with $\Ov{\vr} = \infty$. In the case 
of the hard--sphere EOS, we even obtain a positive lower bound on $\tvu$.

The second observation is a comparison principle. Suppose 
$\tvu_1$, $\tvu_2$ are two solutions of the initial value problem 
\eqref{s17} with the same initial datum but different constants 
\[
\Lambda_1 > \Lambda_2. 
\]
Then 
\begin{equation} \label{s20}
	\tvu_1 > \tvu_2 \ \mbox{in}\ (0,1].
\end{equation}
Indeed, as $\Lambda_1 > \Lambda_2$, we have $\tvu_1 > \tvu_2$ on 
{$(0, \delta)$} for some $\delta > 0$. Let $\delta$ by maximal with this property, meaning 
\[
\tvu_1 > \tvu_2 \ {\mbox{in}\ (0, \delta)},\ 
\tvu_1 (\delta) = \tvu_2 (\delta). 
\]
Necessarily, 
\[
{\partial_r (\tvu_1 - \tvu_2)(\delta) \leq 0}, 
\]
in contrast with equation \eqref{s17} yielding 
\[
\partial_r (\tvu_1 - \tvu_2) (\delta) = \left( \frac{4}{3} \mu + {\lambda} \right)^{-1}(\Lambda_1 - \Lambda_2) > 0. 
\]
We conclude there is at most one solution satisfying the boundary condition
\[
\tvu (1) = u^+_B.
\]

Choosing 
\[
\Lambda + p (\vr_B) + \vr_b (u_B^-)^2 = 0
\]
we obtain the constant solution $\tvu = u^-_B$. It is easy to check that one can solve the boundary value problem for any 
$u_B^+ > u_B^-$ choosing the integration constant $\Lambda$ large enough. Moreover, as $\partial_r \tvu$ satisfies a homogeneous first order linear equation, we get 
\begin{equation} \label{s21}
	\partial_r \tvu > 0.
\end{equation}	
Finally, going back to \eqref{s15} we conclude
\begin{equation} \label{s22}
	\partial_r \tvr = - \frac{\tvr}{\tvu} \partial_r \tvu 
	\approx - \frac{\vr_B}{u^-_B} \partial_r \tvu
\end{equation}
provided 
\begin{equation} \label{s23}
	0 < u^+_B - u^-_B \leq \delta, 
\end{equation}
$\delta > 0$ small enough.

\subsubsection{Inflow--outflow problem in the radially symmetric setting}
\label{ss2}

Finally, we consider radially symmetric solutions of the 
system \eqref{S1}, \eqref{S2} on the spatial domain \eqref{p8}, with the
inflow--outflow boundary conditions \eqref{p9}, \eqref{p10}. Specifically, we set 
\[
\vr_s(x) = \tvr(|x|) \equiv \tvr(r),\ 
\vu_s(x) = \frac{x}{|x|} \tvu (|x|),\ \tvu(|x|) = \tvu(r),\ 
r^- \leq r \leq r^+, 
\]
where 
\begin{align} 
	\partial_r (r^2 \tvr \tvu) &= 0, \label{s6}\\
	\tvr \tvu \partial_r \tvu + \partial_r p(\tvr) &= 
	\left( \frac{4}{3} \mu + {\lambda} \right) \partial_r \left( \frac{\partial_r (r^2 \tvu)  }{r^2}     \right)
	\label{s7}, 
\end{align}
supplemented with the boundary conditions
\begin{align}
	\tvu (r^-) = u^-_B,\ \tvr (r^-) = \vr_B,\label{s8}\\
	\tvu (r^+) = u^+_B.
\label{s9}
\end{align}	 

It follows from the equation of continuity \eqref{s6}, combined with the boundary conditions \eqref{s8}, 
\begin{equation} \label{s10} 
r^2 \tvr \tvu = (r^-)^2 u_B^- \vr_B. 	
\end{equation}
In particular, the momentum equation \eqref{s7} reduces to 
\begin{equation} \label{s11}
	\left( \frac{4}{3} \mu + \lambda \right) 
	\left( \partial^2_{r,r} \tvu + 2 \partial_r \left( \frac{\tvu}{r} \right) \right) = \partial_r p \left( \frac{ (r^-)^2 u_B^- \vr_B }{r^2 \tvu}  \right) + \frac{ (r^-)^2 u_B^- \vr_B }{r^2} \partial_r \tvu.
\end{equation}

\subsection{Main results}
\label{mm}

Having collected the preparatory material we are ready to state our main results. The first concerns the exterior domain outflow problem.

\begin{Theorem}[{\bf Exterior outflow problem}] \label{Tm1}
	Let the pressure satisfy \eqref{p3}, \eqref{m7}, with $0 < \Ov{\vr} < \infty$ (the hard--sphere EOS). Let $\vr_s$, $\vu_s$ be the radially symmetric solution of the 
	stationary problem \eqref{S1}, \eqref{S2} on the exterior domain \eqref{p5} satisfying the outflow boundary conditions \eqref{p6} and the far field conditions \eqref{p7}, where $0 < \vr_\infty < \Ov{\vr}$.
	
	Then there exists $\delta > 0$ such that any weak solution $(\vr, \vu)$ of the corresponding evolutionary NS system converges to $(\vr_s, \vu_s)$ as $t \to \infty$, specifically, 
	\begin{align}
	\vr(t, \cdot) \to \vr_s \ \mbox{in}\ L^q_{\rm loc}(\Omega),\ 1 \leq q < \infty, \br
	\vr \vu (t, \cdot) \to \vr_s \vu_s \ \mbox{weakly in}\ L^2(\Omega; R^3)	
		\label{s38}
	\end{align}	 
as $t \to \infty$ as long as 
\[
0 < u_B	< \delta.
\]
	\end{Theorem}
	
It is worth noting that the convergence is unconditional, meaning it holds for any finite energy weak solution that need not be radially symmetric.

The second result concerns the inflow--outflow on the periodic strip. 

 	\begin{Theorem}[{\bf Inflow--outflow problem on a periodic strip}] \label{Tm2}
 		Let the pressure satisfy \eqref{p3}, \eqref{m7}, and, in addition, {$\Ov{\vr} \leq \infty$},
 		\begin{equation} \label{s39}
 		\liminf_{\vr \to \Ov{\vr}-} \frac{p(\vr)}{P(\vr) + \vr^\gamma} > 0,\ 
 		\gamma > \frac{d}{2}
 		\end{equation}
Let the domain $\Omega$ be the periodic strip specified in \eqref{s12}. Suppose the (constant) inflow density $0 < \vr_B < \Ov{\vr}$ and the velocity $u_B^-$ are given.  		
 		
Let $(\vr_s, \vu_s)$ be the solution of the 
 		stationary problem \eqref{S1}, \eqref{S2} in $\Omega$ satisfying the 
 		boundary conditions \eqref{s13}, \eqref{s14}.
Let  $(\vr, \vu)$ be a weak solution of the corresponding evolutionary NS system satisfying the differential form of the energy inequality \eqref{m6}. 		
 		
Then there exists $\delta > 0$ such that
 		\begin{align}
 			\intO{ \Big[ P(\vr)(t, \cdot) - P'(\vr_s)(\vr(t, \cdot) - \vr_s) - P(\vr_s) \Big] } \to 0, \br 
 			\intO{ \vr(t, \cdot) |\vu (t,\cdot) - \vu_s|^2 } \to 0 
 			\label{s40}
 		\end{align}	 
 		as $t \to \infty$ as long as 
 		\[
 		0 < u^+_B - u^-_B	< \delta.
 		\]
 	\end{Theorem}
 	
Note that \eqref{s39} is satisfied for the isentropic pressure 
$p(\vr) = a \vr^\gamma$, $\gamma > \frac{d}{2}$, where $\Ov{\vr} = \infty$. The convergence stated in \eqref{s40} implies
\begin{align} 
\vr(t, \cdot) &\to \vr_s \ \mbox{in}\ L^\gamma (\Omega),\ \br 
\vr \vu (t, \cdot) &\to \vr_s \vu_s \ \mbox{in}\ L^{\frac{2 \gamma}{\gamma + 1}}(\Omega; R^3)	
	\label{s41}
\end{align} 	
as $t \to \infty$. 	

The last results concerns the inflow--outflow problem on the annular domain.
 
 	\begin{Theorem}[{\bf Inflow--outflow problem on annular domains}] \label{Tm3}
	Let the pressure satisfy \eqref{p3}, \eqref{m7}, and, in addition, {$\Ov{\vr} \leq \infty$}, 
	\begin{equation} \label{s42}
	{	\liminf_{\vr \to \Ov{\vr}-} \frac{p(\vr)}{P(\vr) + \vr^\gamma} > 0. }
	\end{equation}
	Let $\Omega$ be the annular domain \eqref{p8}, where 
	\[
	0 < r^- < r^+ = r^- + 1.
	\]
	Suppose the (constant) inflow density $0< \vr_B < \Ov{\vr}$ and the velocity $u_B^-$ are given.  		
	
	Let $(\vr_s, \vu_s)$ be the radially symmetric solution of the 
	stationary problem \eqref{S1}, \eqref{S2} in $\Omega$ satisfying the 
	boundary conditions \eqref{p9}, \eqref{p10}.
	
	Let  $(\vr, \vu)$ be a weak solution of the corresponding evolutionary NS system satisfying the differential form of the energy inequality \eqref{m6}. 		
	
	Then there exists $\delta > 0$ such that
		\begin{align}
		\intO{ \Big[ P(\vr)(t, \cdot) - P'(\vr_s)(\vr(t, \cdot) - \vr_s) - P(\vr_s) \Big] } \to 0, \br 
		\intO{ \vr(t, \cdot) |\vu (t,\cdot) - \vu_s|^2 } \to 0 
		\label{s43}
	\end{align}	 
	as $t \to \infty$ as long as 
	\[
	0 < u^+_B - u^-_B	< \delta,\ r^- \geq r(\delta).
	\]
\end{Theorem}

{The rest of the paper is devoted to the proofs of the above results.}

\section{Relative energy and relative energy inequality}		
\label{r}

We recall the relative energy inequality satisfied by any weak solution of 
NS system specified in Section \ref{wf}. 

\subsection{Relative energy}

The relative energy is defined as 
\begin{equation} \label{r1}
	E \left( \vr, \vu \Big| R , \vU \right) = 
\frac{1}{2} \vr |\vu - \vU|^2 + P(\vr) - P'(R)(\vr - R) - 
P(R).
\end{equation}	
As $p$ is strictly increasing, the pressure potential $P$ is strictly convex, 
{the energy 
\[
E(\vr, \vm) = \frac{1}{2} \frac{|\vm|^2}{\vr} + P(\vr) 
\]
is a strictly convex function of $(\vr, \vm)$; }
whence the relative energy represents a Bregman divergence between the pairs
$(m = \vr \vu,r)$ and $(M = r\vU, R)$, see \cite[Chapter 3, Section 3.1]{FeiNovOpen}. {Unlike in the energy inequality \eqref{m5}, the quantities $M$ and $R$ are arbitrary smooth functions of 
$t$ and $x$.}	

\subsection{Relative energy inequality}

The weak solutions of the NS system satisfying the integrated variant of the 
energy inequality \eqref{m5} satisfy the {\emph{relative energy inequality}} in the form:
\begin{align} 
	&\left[ \intO{ E \left(\vr, \vu \Big| R, \vU \right) } \right]_{t = 0}^{t = \tau} 
	+ \int_0^\tau {\intO{ \mathbb{S} \Big(\Ds\vu - \Ds \vU \Big): \Big( \Ds \vu - \Ds \vU \Big) } \dt} \br 
	&+ \int_{\partial \Omega} \Big( P(\vr) - P'(R) (\vr - R) - P(R) \Big) [\vuB \cdot \vc{n}]^+ \D \sigma_x \dt 	\br 
	&+ \int_0^\tau \int_{\partial \Omega} \Big( P(\vr_B) - P'(R) (\vr_B - R) - P(R) \Big) [\vuB \cdot \vc{n}]^- \D \sigma_x \dt \br
	&\leq - \int_0^\tau \intO{ \vr (\vu - \vU) \cdot \Grad \vU \cdot (\vu - \vU) } \dt
	-   \int_0^\tau \intO{ \Big( p(\vr) - p'(R) (\vr - R) - p(R) \Big) \Div \vU } \dt \br 
	&\quad + \int_0^\tau \intO{ \frac{\vr}{R} (\vU - \vu) \cdot \Big[ \partial_t (R \vU) + \Div (R \vU \otimes \vU) + \Grad p(R) - {\Div \mathbb{S}(\Ds \vU)}   \Big] } \dt \br 
	&\quad + \int_0^\tau \intO{ \Big( \frac{\vr}{R} (\vu - \vU) \cdot \vU + 
		p'(R) \left( 1 - \frac{\vr}{R} \right) \Big) \Big[ \partial_t R + \Div (R \vU) \Big] } \dt \br 
	&\quad + \int_0^\tau \intO{ \left( \frac{\vr}{R} - 1 \right) (\vU - \vu) \cdot 
		{\Div \mathbb{S} (\Ds \vU) } } \dt	
	\label{r3i}	
\end{align}
for any $\tau \geq 0$. Here $(\vr, \vu)$ is a weak solution while 
$(R,U)$ is a pair of sufficiently smooth ``test'' functions, 
\begin{equation} \label{r3b}
	R > 0,\ {\vc{U}|_{\partial \Omega} = \vuB|_{\partial \Omega}}. 
\end{equation}	
In the case of an exterior domain, we also require the test functions to 
approach the far field state sufficiently fast 
\begin{equation} \label{r3c}
	\vU \to 0,\ R \to \vr_\infty \ \mbox{as}\ |x| \to \infty
\end{equation}
so that all relevant integrals are finite. {Note carefully that the energy inequality \eqref{m5}, that is an integral part of the 
\emph{definition} of weak solution, is a special case of the relative energy inequality \eqref{r3i} for $r = \vr_\infty$, $\vc{U} = \vu_\infty$. 
Still \eqref{r3i} can be derived from \eqref{m5} and the two field equations \eqref{m1}, \eqref{m3}, see e.g. 
\cite[Chapter 3]{FeiNovOpen}.
}		

{If a weak solution satisfies the differential form of the energy inequality \eqref{m6}, the 
	differential form of the
relative energy inequality holds:}
\begin{align} 
\frac{\D}{\dt} &\intO{ E \left(\vr, \vu \Big| R, \vU \right) } 
+ \intO{ \mathbb{S} {\Big(\Ds \vu - \Ds \vU \Big): \Big( \Ds \vu - \Ds \vU \Big) } } \br 
&+ {\int_{\partial \Omega} \Big( P(\vr) - P'(R) (\vr - R) - P(R) \Big) [\vuB \cdot \vc{n}]^+ \D \sigma_x} 	\br 
&+ {\int_{\partial \Omega} \Big( P(\vr_B) - P'(R) (\vr_B - R) - P(R) \Big) [\vuB \cdot \vc{n}]^-  \D \sigma_x} \br
&\leq - \intO{ \vr (\vu - \vU) \cdot \Grad \vU \cdot (\vu - \vU) }
-   \intO{ \Big( p(\vr) - p'(R) (\vr - R) - p(R) \Big) \Div \vU } \br 
&+ \intO{ \frac{\vr}{R} (\vU - \vu) \cdot \Big[ \partial_t (R \vU) + \Div (R \vU \otimes \vU) + \Grad p(R) - {\Div \mathbb{S}(\Ds \vU) }   \Big] } \br 
&+ \intO{ \Big( \frac{\vr}{R} (\vu - \vU) \cdot \vU + 
	p'(R) \left( 1 - \frac{\vr}{R} \right) \Big) \Big[ \partial_t R + \Div (R \vU) \Big] } \br 
&+ \intO{ \left( \frac{\vr}{R} - 1 \right) (\vU - \vu) \cdot 
	{\Div \mathbb{S} (\Ds \vU) } }
\label{r3}	
	\end{align}
{in $\mathcal{D}'(0, \infty)$.}

\section{Convergence via relative energy}
\label{c}

A natural idea is to show the convergence results claimed in Section \ref{mm}
by using the steady state solution as a test function in the relative energy inequality.

\subsection{Proof of Theorem \ref{m1}}

The relative energy inequality \eqref{r3i} evaluated for the pairs 
$(\vr, \vu)$, $R = \vr_s$, $\vU = \vu_s$ gives rise to 
\begin{align} 
	&\left[ \intO{ E \left(\vr, \vu \Big| \vr_s, \vu_s \right) } \right]_{t = 0}^{t = \tau} +  
	+ \int_0^\tau \intO{ \mathbb{S} {\Big(\Ds \vu - \Ds \vu_s \Big): \Big( \Ds \vu - \Ds \vu_s \Big) } } \dt \br 
	&+ u_B \int_0^\tau \int_{|x|=r} \Big( P(\vr) - P'(\vr_s) (\vr - \vr_s) - P(\vr_s) \Big)  \D \sigma_x \dt 	\br 
	&\leq - \int_0^\tau \intO{ \vr (\vu - \vu_s) \cdot \Grad \vu_s \cdot (\vu - \vu_s) } \dt\br
	&\quad -   \int_0^\tau \intO{ \Big( p(\vr) - p'(\vr_s) (\vr - \vr_s) - p(\vr_s) \Big) \Div \vu_s } \dt \br 
		&\quad + \int_0^\tau \intO{ \left( \frac{\vr}{\vr_s} - 1 \right) (\vu_s - \vu) \cdot 
	\Big( \vr_s \vu_s \cdot \Grad \vu_s + \Grad p(\vr_s) \Big) } \dt	.
	\label{c1}	
\end{align}

Next, as the pressure satisfies the hard--sphere EOS, we have 
\begin{equation} \label{c2}
0 \leq \vr (t,x) \leq \Ov{\vr} { < \infty} 
\end{equation}
and, in accordance with \eqref{s2},
\begin{equation} \label{c3}
 \left( \frac{\vr}{\vr_s} - 1 \right)^2 \aleq 
\Big( p(\vr) - p'(\vr_s) (\vr - \vr_s) - p(\vr_s) \Big)  
\end{equation}
as long as $u_B> 0$ is small enough.

{Seeing that
\[	  
(\vu - \vu_s)|_{\partial \Omega} = 0
\]
we have Korn's inequality 
\[
\intO{ |\Grad \vu - \Grad \vu_s|^2 } \aleq \intO{ \mathbb{S} \Big(\Ds \vu - \Ds \vu_s \Big): \Big( \Ds \vu - \Ds \vu_s \Big) }.
\]
}	 
Thus, by Hardy's inequality, 
\begin{equation} \label{c4}
\intO{ \frac{|\vu - \vu_s|^2}{|x|^2} } \aleq \intO{ {\mathbb{S} \Big(\Ds \vu - \Ds \vu_s \Big): \Big( \Ds \vu - \Ds \vu_s \Big) } }.
	\end{equation}
Consequently, in view of the uniform bounds \eqref{s3}
\begin{equation} \label{c5}
| \Grad \vu_s | \aleq u_B \frac{1}{|x|^3},
	\end{equation}	
and the integral 
\[
\left| \intO{ \vr (\vu - \vu_s) \cdot \Grad \vu_s \cdot (\vu - \vu_s) } \right| \aleq \Ov{\vr} u_B \intO{ \frac{1}{|x|^2} |\vu - \vu_s|^2 }
\]	
can be absorbed by the dissipative term 
\[
\intO{ \mathbb{S} {\Big(\Ds \vu - \Ds \vu_s \Big): \Big( \Ds \vu - \Ds \vu_s \Big) } } 
\]
as long as $u_B$ is small enough. Accordingly, inequality \eqref{c1} reduces to 
\begin{align} 
	&\left[ \intO{ E \left(\vr, \vu \Big| \vr_s, \vu_s \right) } \right]_{t = 0}^{t = \tau} 
	+ \frac{1}{2} \int_0^\tau \intO{ \mathbb{S} {\Big(\Ds \vu - \Ds \vu_s \Big): \Big( \Ds \vu - \Ds \vu_s \Big) } } \dt \br 
	&+ u_B \int_0^\tau \int_{|x|=r} \Big( P(\vr) - P'(\vr_s) (\vr - \vr_s) - P(\vr_s) \Big) \D \sigma_x \dt 	\br 
	&
	+  \int_0^\tau \intO{ \Big( p(\vr) - p'(\vr_s) (\vr - \vr_s) - p(\vr_s) \Big) \Div \vu_s } \dt \br 
	&\leq \int_0^\tau \intO{ \left( \frac{\vr}{\vr_s} - 1 \right) (\vu_s - \vu) \cdot 
		\Big( \vr_s \vu_s \cdot \Grad \vu_s + \Grad p(\vr_s) \Big) }	.
	\label{c6}	
\end{align}

Now, going back to the relation \eqref{s5}, {and using \eqref{s2}, \eqref{s3}}, we deduce 
\begin{equation} \label{c7}
	\Div \vu_B \ageq \frac{1}{\vr_\infty} u_b^3 \frac{1}{|x|^7},  
\end{equation}
while 
\begin{equation} \label{c8}
\Big| \vr_s \vu_s \cdot \Grad \vu_s + \Grad p(\vr_s) \Big| \aleq 
u^2_B \frac{1}{|x|^5}.
\end{equation} 
Consequently, by virtue of \eqref{c3}, \eqref{c7},
\begin{align} 
&\left| \intO{ \left( \frac{\vr}{\vr_s} - 1 \right) (\vu_s - \vu) \cdot 
	\Big( \vr_s \vu_s \cdot \Grad \vu_s + \Grad p(\vr_s) \Big) }
\right| \br 
&\quad \aleq \left| \intO{ \left( \frac{\vr}{\vr_s} - 1 \right) \frac{|\vu_s - \vu|}{|x|} \cdot u^2_B \frac{1}{|x|^4}
 }
\right| \br &\quad \aleq 
\intO{ \left( \frac{\vr}{\vr_s} - 1 \right)^2 u_B^{\frac{7}{2}} \frac{1}{|x|^8} } + \intO{ \frac{|\vu_s - \vu|^2}{|x|^2} u_b^{\frac{1}{2}} } \br 
&\quad \aleq 
u_B^{\frac{1}{2}}  \intO{ \Big( p(\vr) - p'(\vr_s) (\vr - \vr_s) - p(\vr_s) \Big) \Div \vu_s } \br &\quad + u_B^{\frac{1}{2}} {\intO{ \mathbb{S} \Big(\Ds \vu - \Ds \vu_s \Big): \Big( \Ds \vu - \Ds \vu_s \Big) }}
\nonumber
\end{align}
Thus we conclude 
\begin{align} 
	&\left[ \intO{ E \left(\vr, \vu \Big| \vr_s, \vu_s \right) } \right]_{t = 0}^{t = \tau} 
	+ \frac{1}{4} \int_0^\tau \intO{ \mathbb{S} {\Big(\Ds \vu - \Ds \vu_s \Big): \Big( \Ds \vu - \Ds \vu_s \Big) } } \dt \br 
	&+ u_B \int_0^\tau \int_{|x|=r} \Big( P(\vr) - P'(\vr_s) (\vr - \vr_s) - P(\vr_s) \Big) \D \sigma_x \dt 	\br 
	&
	+  \frac{u_B^3}{\vr_\infty} \int_0^\tau \intO{ \Big( p(\vr) - p'(\vr_s) (\vr - \vr_s) - p(\vr_s) \Big) \frac{1}{|x|^7}  } \dt \leq 0
	\label{c9}	
\end{align}
provided $0 < u_B < \delta$,  $\delta > 0$ small enough.

It follows 
\begin{equation} \label{c9a}
\sup_{\tau \geq 0} \intO{ E \left( \vr, \vu \Big| \vr_s, \vu_s \right) } \aleq 1,
\end{equation}
\begin{equation} \label{c9b}
\int_0^\infty \intO{ {|\Grad \vu - \Grad \vu_s |^2 } } \dt \aleq 1,
\end{equation}		
and, finally, 
\begin{equation} \label{c9c}
\int_0^\infty 	\intO{ \Big( p(\vr) - p'(\vr_s) (\vr - \vr_s) - p(\vr_s) \Big) \frac{1}{|x|^7}  } \dt \aleq 1.
\end{equation}
Moreover, since the pressure is strictly convex and the density bounded above by $\Ov{\vr}$, the bound \eqref{c9c} yields 
\begin{equation} \label{c9d}
	\int_0^\infty 	\intO{ |\vr - \vr_s |^2 \frac{1}{|x|^7}  } \dt \aleq 1.
\end{equation}

{As shown in \cite[Proposition 2.2]{FanFeiHof}, we have 
	\begin{equation} \label{c9e}
		\vr \in C([\tau_1, \tau_2]; L^1_{\rm loc}({\Omega}))
	\end{equation}
	for any $0 < \tau_1 < \tau_2$. In particular, as $\vr \leq \Ov{\vr}$, we deduce 
	\begin{equation} \label{c9f}
		b(\vr) \in C([\tau_1, \tau_2]; L^1_{\rm loc}({\Omega}))	
	\end{equation}
	for any $b \in C(R)$. Note carefully that the weak and renormalized formulation of the 
equation of continuity \eqref{m1}, \eqref{m2} imply the existence of a (weakly in times) continuous  
representative of $\vr$ and $\widetilde{b(\vr)}$, respectively. Relations \eqref{c9e}, \eqref{c9f} 
rigorously justify the seemingly obvious conclusion
\[
b(\vr)(t, \cdot) = \widetilde{b(\vr)}(t, \cdot) \ \mbox{for \emph{any}}\ t \geq 0.
\]
}

{Next, using the weak formulation of the continuity equation \eqref{m1}, and estimates \eqref{c9b}, \eqref{c4}, together with the bound $0<\vr\le \Ov{\vr}$, we conclude that the sequence of functions 
\begin{equation}\label{uni-bound}
\left\{t\mapsto \intO{ \vr(n+t) \psi } \right\}_{n\in N} \ \mbox{is  uniformly equi--continuous  for each}\ \psi\in C^\infty_c(\Omega), t\in[0,1].
\end{equation}
}

Since $C^\infty_c(\Omega)$ is dense in {$L^{q'}(\Omega)$, $1 < q<\infty,$ $1/q+1/q'=1$  and $\{\vr(n+\cdot)\}$ are uniformly bounded, it follows that 
\begin{equation}\label{con-vr}
\intO{ \vr_n \psi } \to \vr_\psi^\infty \ \mbox{uniformly for}\ t \in [0,1] \ \mbox{for any compactly supported}\ 
\psi \in L^{q'}(\Omega)
\end{equation} 
passing to a subsequence as the case may be. }
In addition,  by \eqref{c9d},  
\begin{equation} \label{c9ee}
	\int_n^{(n+1)} 	\intO{ |\vr - \vr_s |^2 \frac{1}{|x|^7}  } \dt \to 0
\end{equation}
as $n\to\infty$. Consequently,
\begin{equation}
\vr^\infty_\psi(t)= \intO{ \vr_s\psi },
\end{equation}
in particular, the convergence is unconditional. Finally, using now the renormalized version of the continuity equation 
\eqref{m2}
in the same fashion we obtain 
\begin{equation}
\intO{ |\vr(n+t)|^q \psi } \to \intO{ |\vr_s|^q\psi } \ \mbox{in}\ C[0,1], 1 \leq q < \infty, 
\end{equation}
{all compactly supported $\psi \in L^1(\Omega)$.}
Thus, with help of  Radon-Riesz theorem, we conclude from weak convergence and convergence of norms that
\begin{equation}
\lim_{t\to\infty}\|\vr(t,\cdot)-\vr_s(\cdot)\|_{L^q_{\rm loc}(\Omega)}=0.
\end{equation}

The second part of \eqref{s38} claiming the weak convergence of momenta can be shown in a similar fashion using the weak continuity of $\vr \vu$. Indeed it follows from the weak formulation of momentum equation \eqref{m3} that 
\begin{equation}\label{uni-bound-momentum}
\left\{ t\mapsto\intO{ \vm(n+t) \cdot \psi } \right\}_{n\in N} \mbox{ is  uniformly equi--continuous  for each } \psi\in C^\infty_c(\Omega;R^3), t\in[0,1]
\end{equation}
where $\vm=\vr\vu$. Moreover, from the bounds $0\le\vr(t,x)\le \Ov{\vr}$ and {$\vr|\vu - \vu_s|^2\in L^\infty(0,\infty;L^1(\Omega))$} we conclude that $\vr(\vu - \vr_s)\in L^\infty(0,\infty;L^2(\Omega))$. {The proof of (weak) convergence claimed 
in \eqref{s38} is then concluded exactly as in the above replacing \eqref{c9ee} by \eqref{c9b}.}
	   
We have proved Theorem \ref{Tm1}. 

\subsection{Proof of Theorem \ref{Tm2}}
\label{sc3}

As the weak solution satisfies the differential form of the energy inequality,
we may write
the relative energy inequality evaluated in terms of the steady solution $\vr_s$, $\vu_s$ as 
\begin{align} 
	\frac{\D}{\dt} &\intO{ E \left(\vr, \vu \Big| \vr_s, \vu_s \right) } 
	+ \intO{ {\mathbb{S} \Big(\Ds \vu - \Ds \vu_s \Big): \Big( \Ds \vu - \Ds \vu_s \Big) }  } \br 
	&+ u^+_B\int_{r = 1} \Big( P(\vr) - P'(\vr_s) (\vr - \vr_s) - P(\vr_s) \Big)  \D \sigma_x 	\br 
	&\leq - \intO{ \vr [u^3 - u^3_s]^2 \partial_r \tvu  }
	-   \intO{ \Big( p(\vr) - p'(\vr_s) (\vr - \vr_s) - p(\vr_s) \Big) \partial_r \tvu } \br 
	&+ \intO{ \left( \frac{\vr}{\vr_s} - 1 \right) (u_s^3 - u^3)  
		\Big( \vr_s \tvu \partial_r \tvu + \partial_r p(\vr_s) \Big) }	.
	\label{ccc1}	
\end{align}

Going back to the estimates \eqref{s21}--\eqref{s23} we may rewrite the above inequality in the form
\begin{align} 
	\frac{\D}{\dt} &\intO{ E \left(\vr, \vu \Big| \vr_s, \vu_s \right) } 
	+ \intO{ {\mathbb{S} \Big(\Ds \vu - \Ds \vu_s \Big): \Big( \Ds \vu - \Ds \vu_s \Big) } } \br 
	&+ \int_{z= 1} \Big( P(\vr) - P'(\vr_s) (\vr - \vr_s) - P(\vr_s) \Big) u_B^+ \D \sigma_x 	\br 
	&+ \intO{ \vr [\vu^3 - \vu^3_s]^2 \partial_z \tvu  }
	+  \intO{ \Big( p(\vr) - p'(\vr_s) (\vr - \vr_s) - p(\vr_s) \Big) \partial_r \tvu } \br 
	&\leq C \intO{ | \vr - \vr_s | |\vu_s^3 - \vu^3 |  
		\partial_r \tvu }	
	\label{ccc2}	
\end{align}
provided $\vr_B$, $u^-_B$ are fixed and $0 < u^+_B - u^-_B < \delta$, where $\delta > 0$ is sufficiently small. 

At this stage, the proof can be completed in the same way as that of Theorem 
\ref{Tm1} in the preceding section as long as the pressure satisfies the 
hard--sphere EOS. {We only remark that thanks to the differential form of the relative energy inequality, 
the relative energy 
\[
\intO{ E \left(\vr, \vu \Big| \vr_s, \vu_s \right) (t, \cdot) } 
\]
is a non--increasing function of time, in particular, it admits a limit $\mathcal{E}_\infty \geq 0$ as $t \to \infty$. Moreover, it follows from 
\eqref{c9b}, \eqref{c9c} that, necessarily, $\mathcal{E}_\infty = 0$, which yields the strong convergence claimed in 
\eqref{s40}.  }	 

Now, we focus on the case $\Ov{\vr} = \infty$ compatible with the isentropic pressure law $p = a \vr^\gamma$, $\gamma > \frac{d}{2}$.
First, we estimate the integral on the right--hand side for large values of density, $\vr \geq \Ov{\vr}_s > \vr_s$: 
\begin{align}
	C &\int_{\vr > \Ov{\vr}_s} | \vr - \vr_s | |u_s^3 - u^3 |  
	\partial_r \tvu \ \dx \leq C \int_{\vr > \Ov{\vr}_s} \vr  |u_s^3 - u^3 |  
	\partial_r \tvu \ \dx \br &\leq 
	\frac{1}{2} \intO{ \vr [u^3 - u^3_s]^2 \partial_r \tvu  } + 
	C_2 \int_{\vr \geq \Ov{\vr}_s} \vr \partial_r \tvu \dx \br 
	&\leq 
	\frac{1}{2} \intO{ \vr [ u^3 - u^3_s]^2 \partial_r \tvu  } \br&+ 
	C_2 \int_{\vr \geq \Ov{\vr}_s} \frac{\vr}{  p(\vr) - p'(\vr_s) (\vr - \vr_s) - p(\vr_s)} \Big( p(\vr) - p'(\vr_s) (\vr - \vr_s) - p(\vr_s) \Big) \partial_r \tvu \dx,
	\nonumber
\end{align}
where $C_2$ depends only on $C_1$ but not on $\Ov{\vr}_s$. 
In view of hypothesis \eqref{s39}, 
\[
	\liminf_{\vr \to \infty} \frac{p(\vr)}{\vr} = \infty,
\]
Thus we may choose $\Ov{\vr}_s > 0$ large enough so that 
\begin{align}
	C_2 &\int_{\vr \geq \Ov{\vr}_s} \frac{\vr}{  p(\vr) - p'(\vr_s) (\vr - \vr_s) - p(\vr_s)} \Big( p(\vr) - p'(\vr_s) (\vr - \vr_s) - p(\vr_s) \Big) \partial_r \tvu \dx \br
	&\leq \frac{1}{2}\intO{ \Big( p(\vr) - p'(\vr_s) (\vr - \vr_s) - p(\vr_s) \Big) \partial_r \tvu }.
	\nonumber
\end{align}
Going back to \eqref{ccc2} we therefore conclude 
\begin{align} 
	\frac{\D}{\dt} &\intO{ E \left(\vr, \vu \Big| \vr_s, \vu_s \right) } 
	+ \intO{ \mathbb{S} {\Big(\Ds \vu - \Ds \vu_s \Big): \Big( \Ds \vu - \Ds \vu_s \Big) } } \br 
	&+ u^+_B \int_{r = 1} \Big( P(\vr) - P'(\vr_s) (\vr - \vr_s) - P(\vr_s) \Big) \D \sigma_x 	\br 
	&+ \frac{1}{2} \intO{ \vr [u^3 - u^3_s]^2 \partial_r \tvu  }
	+  \frac{1}{2} \intO{ \Big( p(\vr) - p'(\vr_s) (\vr - \vr_s) - p(\vr_s) \Big) \partial_r \tvu } \br 
	&\leq C \int_{\vr \leq \Ov{\vr}_s} | \vr - \vr_s | |u_s^3 - u^3 |  
	\partial_r \tvu \ \dx
	\label{ccc4}	
\end{align}
for a certain $\Ov{\vr}_s > 0$ large enough. 

As the pressure is strictly convex, we get 
\begin{equation} \label{ccc5}
	|\vr - \vr_s|^2 \leq c(\Ov{\vr}_s) \Big( p(\vr) - p'(\vr_s) (\vr - \vr_s) - p(\vr_s) \Big) 
\end{equation} 
whenever $\vr, \vr_s \leq \Ov{\vr}_s $. Thus making use of Poincar\' e inequality we may control the last integral in \eqref{ccc4} 
as
\begin{align}
	C &\int_{\vr \leq \Ov{\vr}_s} | \vr - \vr_s | |\vu_s^3 - \vu^3 |  
	\partial_r \tvu \ \dx \leq \frac{C}{2} \int_{\vr \leq \Ov{\vr}_s} 
	|\vr - \vr_s|^2 (\partial_r \tvu)^{\frac{3}{2}} + \frac{C}{2} 
	\intO{ | \vu - \vu_s |^2 	(\partial_r \tvu)^{\frac{1}{2}}}\br &\leq 
	C c(\Ov{\vr}_s) \intO{\Big( p(\vr) - p'(\vr_s) (\vr - \vr_s) - p(\vr_s) \Big) (\partial_z \tvu)^{\frac{3}{2}} } \br
	&+ \frac{C}{2} c_P \| \partial_z \tvu \|^{\frac{1}{2}}_{L^\infty}
	\intO{ \mathbb{S} \Big(\Grad \vu - \Grad \vu_s \Big): \Big( \Grad \vu - \Grad \vu_s \Big) }
	\nonumber
\end{align}	
Thus the right--hand side of \eqref{ccc4} may be absorbed by the left--hand side as long as $\delta > 0$ is small enough. 
We conclude 
\begin{align} 
	\frac{\D}{\dt} &\intO{ E \left(\vr, \vu \Big| \vr_s, \vu_s \right) } 
	+ \intO{ {\mathbb{S} \Big(\Ds \vu - \Ds \vu_s \Big): \Big( \Ds \vu - \Ds \vu_s \Big) } } \br 
	&+ \int_{r = 1} \Big( P(\vr) - P'(\vr_s) (\vr - \vr_s) - P(\vr_s) \Big) u_B^+ \D \sigma_x 	\br 
	&+ \frac{1}{4} \intO{ \vr [u^3 - u^3_s]^2 \partial_r \tvu  }
	+  \frac{1}{4} \intO{ \Big( p(\vr) - p'(\vr_s) (\vr - \vr_s) - p(\vr_s) \Big) \partial_r \tvu } \leq 0.
	\label{ccc6}	
\end{align}

{Now, we complete the proof exactly as in  the case of hard--sphere EOS.}
Relation \eqref{ccc6} implies 
\begin{equation} \label{ccc7}
\intO{ E \left(\vr, \vu \Big| \vr_s, \vu_s \right)(t, \cdot) } \searrow \mathcal{E}_\infty \geq 0 \ \mbox{as}\ t \to \infty, 
\end{equation}
and
\begin{equation} \label{ccc8}
\int_0^\infty \intO{ \left[ \mathbb{S} {\Big(\Ds \vu - \Ds \vu_s \Big): \Big( \Ds \vu - \Ds \vu_s \Big) } + \Big( p(\vr) - p'(\vr_s) (\vr - \vr_s) - p(\vr_s) \Big)\right] } \dt < \infty.
\end{equation}
Consequently, relation \eqref{ccc8}, together with hypothesis \eqref{s39}, 
give rise to
\[ 
\intO{ E \left(\vr, \vu \Big| \vr_s, \vu_s \right)(t_n, \cdot) } \to 0
\]
for some $t_n \to \infty$, which yields the desired conclusion $\mathcal{E}_\infty = 0$ in \eqref{ccc7}.

\subsection{Inflow-outflow problem in the radially symmetric setting}

Our ultimate goal is to prove Theorem \ref{Tm3}. 
Similarly to the preceding section, we have  
\begin{align} 
	\frac{\D}{\dt} &\intO{ E \left(\vr, \vu \Big| \vr_s, \vu_s \right) } 
	+ \intO{ \mathbb{S} {\Big(\Ds \vu - \Ds \vu_s \Big): \Big( \Ds \vu - \Ds \vu_s \Big) } } \br 
	&+ \int_{|x|=r^+} \Big( P(\vr) - P'(\vr_s) (\vr - \vr_s) - P(\vr_s) \Big) u_B^+ \D \sigma_x 	\br 
	&\leq - \intO{ \vr (\vu - \vu_s) \cdot \Grad \vu_s \cdot (\vu - \vu_s) }
	-   \intO{ \Big( p(\vr) - p'(\vr_s) (\vr - \vr_s) - p(\vr_s) \Big) \Div \vu_s } \br 
	&+ \intO{ \left( \frac{\vr}{\vr_s} - 1 \right) (\vu_s - \vu) \cdot 
		\Big( \vr_s \vu_s \cdot \Grad \vu_s + \Grad p(\vr_s) \Big) }	.
	\label{cc1}	
\end{align}
Moreover, formula \eqref{s4} yields 
\begin{equation} \label{cc2}
\xi \cdot \Grad \vu_s \cdot \xi = \partial_r \tvu(|x|) \left|\xi \cdot \frac{x}{|x|} \right|^2+ \frac{\tvu (|x|)}{|x|} \left( |\xi|^2 -  \left|\xi \cdot \frac{x}{|x|} \right|^2 \right).
	\end{equation}
Consequently, the velocity gradient $\Grad \vu_s$ is positively definite as long as 
$\partial_r \tvu > 0$, $\tvu > 0$.

Thus the proof of Theorem \ref{Tm3} will be the same as that of 
Theorem \ref{Tm2} as soon as we observe that the qualitative properties of the 
stationary solutions on the periodic strip and the annular domain are similar. 
{More specifically, this is true as long as the constant $\delta$,} 
\[
0 < u^-_B - u^+_B < \delta
\]
is small and $r = r(\delta)$ is chosen large enough. 

Thus the proof of Theorem \ref{Tm3} reduces to comparison of the solution 
$\tvr_1$, $\tvu_1$ of the ``flat'' problem 
\begin{align}
\partial_s (\tvr_1 \tvu_1) &= 0 , \br 
\tvr_1 \tvu_1 \partial_s \tvu_1 + \partial_s p(\tvr_1) &= \left( \frac{4}{3} \mu + \eta \right) \partial^2_{s,s} \tvu_1 	 
	\label{F1}
\end{align}	
with that of the ``curved'' problem 
\begin{align} 
	\partial_s (s^2 \tvr_2 \tvu_2) &= 0, \br 
\tvr_2 \tvu_2 \partial_s \tvu_2 + \partial_s p(\tvr_2) &= 
\left( \frac{4}{3} \mu + { \lambda} \right) \partial_s \left( \frac{ \partial_s (s^2 \tvu_2)}{s^2} \right).
\label{F2}
\end{align}
Both problems are considered on the interval $s \in [r, r+1]$, with the initial data 
\begin{equation} \label{F3}
\tvr_1(r) = \tvr_2(r) = \vrB, \ \tvu_1(r) = \tvu_2(r) = u^-_B,
\end{equation}
and the boundary data 
\begin{equation} \label{F4}
\tvu_1(r + 1) = \tvu_2(r+1) = u^+_B.
\end{equation}

Going back to Sections \ref{ss3}, \ref{ss2}
we record 
\[
\tvr_1 = \vrB \vu_B^- \frac{1}{\tvu_1},\ \tvr_2 = \frac{r^2}{s^2} \vrB \vuB^- \frac{1}{\tvu_2},\ s \in [r, r+1],
\]
where
\[
\frac{r^2}{s^2} \to 1 \ \mbox{as}\ r \to \infty\  \mbox{uniformly for}
\ s \in [r, r + 1].
\]

Next, the second order term in \eqref{F2} can be written as 
\[
\partial_s \left( \frac{ \partial_s (s^2 \tvu_2)}{s^2} \right) = 
\partial^2_{s,s} \tvu_2 + 2 \partial_s \left( \frac{1}{s} \tvu_2 \right), 
\]
where 
\[
\frac{1}{s} \leq \frac{1}{r} \to 0 \ \mbox{as}\ r \to \infty.
\] 
Consequently, we get 
\[
\tvr_2 \to \tvr_1,\ \tvu_2 \to \tvu_1 \ \mbox{uniformly in}\ 
C^1[r, r + 1] \ \mbox{as}\ r \to \infty.
\]
{The} rest of the proof of Theorem \ref{Tm3} is the same as that of Theorem \ref{Tm2}.

\def\cprime{$'$} \def\ocirc#1{\ifmmode\setbox0=\hbox{$#1$}\dimen0=\ht0
	\advance\dimen0 by1pt\rlap{\hbox to\wd0{\hss\raise\dimen0
			\hbox{\hskip.2em$\scriptscriptstyle\circ$}\hss}}#1\else {\accent"17 #1}\fi}


\end{document}